\documentclass[11pt]{article}
\usepackage{graphicx,amsfonts,amsthm,amssymb, amsmath, epic}
\input xy
\xyoption{all}
\theoremstyle{plain}
\newtheorem{thm}{{\bf Theorem}}[section]
\newtheorem{lemma}{{\bf Lemma}}[section]
\newtheorem{prop}{{\bf Proposition}}[section]

\theoremstyle{remark}
\newtheorem{remark}{{\it Remark}}[section]

\def\be{\begin{eqnarray}}
\def\ee{\end{eqnarray}}
\def\ben{\begin{eqnarray*}}
\def\een{\end{eqnarray*}}
\def\ba{\begin{array}}
\def\ea{\end{array}}
\def\bp{\noindent{\it Proof. }}
\def\ep{\noindent{\hfill \fbox{}}}

\def\definition{\noindent{\bf Definition. }}
\def\pic{{\rm Pic}}

\def\noi{\noindent}

\def\al{\alpha}
\def\ve{\varepsilon}
\def\chal{\alpha^{\vee}}
\def\mc{{\mathbb C}}
\def\mz{{\mathbb Z}}

\def\mpp{{\mathbb P}}
\def\disp{\displaystyle}
\newcommand{\ol}[1]{\overline{#1}}

\newcommand{\wt}[1]{\widetilde{#1}}
\newcommand{\mapright}[1]{%
   \smash{\mathop{%
   \hbox to 1cm{\rightarrowfill}}\limits^{#1}}}
\newcommand{\mapleft}[1]{%
   \smash{\mathop{%
   \hbox to 1cm{\leftarrowfill}}\limits^{#1}}}
\newcommand{\maplleft}[2]{%
   \smash{\mathop{%
   \hbox to 1cm{\leftarrowfill}}\limits_{#1}^{#2}}}

\makeatletter       
\renewcommand{\@biblabel}[1]{#1.}
\makeatother

\begin{document}

\title{Elliptic curves and birational 
representation of Weyl groups}
\author{Eguchi Mitsuaki${}^*$ and Tomoyuki Takenawa${}^{**}$}
\date{}
\maketitle

\begin{center}
${}^*$ Graduate School of Mathematical Sciences, University of Tokyo\\ 
Komaba 3-8-1, Meguro-ku, Tokyo 153-8914, Japan\\

${}^{**}$ Faculty of Marine Technology, Tokyo University of Marine Science
and Technology\\ 
Echujima 2-1-6, Koto-ku, Tokyo 135-8533, Japan\\
\end{center}

\begin{abstract} 
Some Weyl group acts on a family of rational varieties
obtained by successive blow-ups at $m$ ($m\geq n+2$) points in the 
projective space $\mpp^n(\mc)$.
In this paper we study the case where all the points of blow-ups 
lie on a certain elliptic curve in $\mpp^n$.
Investigating the action of Weyl group on the Picard groups 
on the elliptic curve and 
on rational varieties, we show that the action on the parameters 
can be written as a group of linear transformations
on the $(m+1)$-st power of a torus.
\end{abstract} 



\section{Introduction}
\begin{figure}[ht] \label{dynkin1}
\hspace{2cm}
\begin{picture}(300,60)
\put(50,20){\line(1,0){20}}
\put(75,20){\line(1,0){10}}
\dashline{3}(85,20)(110,20)
\put(110,20){\line(1,0){10}}
\put(125,20){\line(1,0){20}}
\put(150,20){\line(1,0){20}}
\put(175,20){\line(1,0){10}}
\dashline{3}(185,20)(210,20)
\put(210,20){\line(1,0){10}}
\put(225,20){\line(1,0){20}}
\put(147.5,22.5){\line(0,1){20}}
\put(47.5,20){\circle{5}}
\put(72.5,20){\circle{5}}
\put(122.5,20){\circle{5}}
\put(147.5,20){\circle{5}}
\put(172.5,20){\circle{5}}
\put(222.5,20){\circle{5}}
\put(247.5,20){\circle{5}}
\put(147.5,45){\circle{5}}
\put(47.5,30){\makebox(0,0){}}
\put(197.5,30){\makebox(0,0){}}
\put(47.5,10){\makebox(0,0){$\al_{1}$}}
\put(72.5,10){\makebox(0,0){$\al_{2}$}}
\put(122.5,10){\makebox(0,0){$\al_{n}$}}
\put(147.5,10){\makebox(0,0){$\al_{n+1}$}}
\put(177.5,10){\makebox(0,0){$\al_{n+2}$}}
\put(227.5,10){\makebox(0,0){$\al_{m-2}$}}
\put(257.5,10){\makebox(0,0){$\al_{m-1}$}}
\put(160,45){\makebox(0,0){$\al_{0}$}}
\end{picture}
\caption[]{$W(n,m)$ Dynkin diagram}
\end{figure}
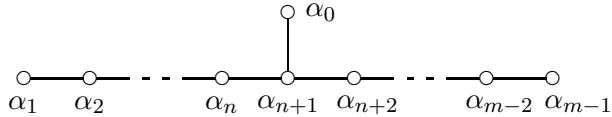

By the works of Coble \cite{coble} and Dolgachev-Ortland \cite{do},
it has been known that
some Weyl group behaves as pseudo-isomorphisms (isomorphisms excluding
sub-varieties of co-dimension $2$ or higher) and acts on
a family of rational varieties
obtained by successive blow-ups at $m$ ($m\geq n+2$) points in $\mpp^n(\mc)$.
Here, the Weyl group is given by the 
Dynkin diagram in Fig. \ref{dynkin1} and denoted by $W(n,m)$. 

Recently, one of the authors (TT)  \cite{takenawa}
introduced dynamical systems defined by
translations of affine Weyl groups (\S 6.5 in \cite{kac}) with the symmetric Cartan matrices
included in $W(n,m)$.
For example, if $m\geq n+7$, $W(n,m)$ includes the affine Weyl group of type 
$E_8^{(1)}$. 
In the case of $n=2, m=9$, this dynamical system coincides with the elliptic
difference Painlev\'e equation proposed by Sakai \cite{sakai},
from which all the discrete and continuous Painlev\'e equations are obtained by
degeneration.
For $n\geq 3$, however, the time evolution of the parameters may not
be solved in general, thus, the systems should be considered not to be 
$n$-dimensional but to be higher dimensional.

On the other hand, Kajiwara {\it et al.} \cite{kmnoy}
has proposed a birational representation of the Weyl group $W(n,m)$,
in which all the points of blow-ups lie on a certain elliptic curve 
in $\mpp^n$. This representation is a special case of \cite{coble,do}.
In this case, 
the action on the parameters can be written as a group of
linear transformations
on the $(m+1)$-st power of a torus. 
The calculation was carried out in a rather heuristic manner in 
\cite{kmnoy}. 

In this paper, we recover the birational representation
of Kajiwara {\it et al.} geometrically,
by investigating the actions of the Weyl group on the Picard groups 
on the elliptic curve and 
on rational varieties. Our method corresponds to
the ``linear map'' or the ``period map'' in 2-dimensional case
\cite{looijenga,sakai}.

This article is organized as follows.
In Section~2, we review the relationship between rational varieties
and groups of Cremona transformations.
In Section~3, it is proved that general two elliptic curves 
in $\mpp^n$ of degree $n+1$
are translated to each other by a projective linear transformation. 
In Section~4, we investigate how the birational representation of the Weyl group
is restricted onto elliptic curves of degree $n+1$ geometrically.
In Section~5, we present some examples of calculation, and 
recover the birational representation
of \cite{kmnoy}. 

\section{The birational representation of Weyl groups by 
Coble and Dolgachev-Ortland}


Let $m\geq n+2 $.  Let $X_{n,m}$ be the 
configuration space of ordered $m$ points in $\mpp^n$:
\begin{equation*} X_{n,m} = 
\mathrm{PGL}(n+1) \left\backslash \left\{
\begin{pmatrix}
a_{01}&\cdots&a_{0m}\\
a_{11}&\cdots&a_{1m}\\
\vdots &\cdots&\vdots \\
a_{n1}&\cdots&a_{nm}\end{pmatrix} \left| 
\ba{l}
\mbox{the determinant of}\\
\mbox{every}\ (n+1)\times(n+1) \\
\mbox{sub-matrix is nonzero}
\ea
\right. \right\}
\right/ (\mc^{\times})^{m} ,
\end{equation*}
which is a quasi-projective variety of dimension $n(m-n-2)$.
We also consider $X_{n,m}^1\simeq X(n,m+1)$ with a natural projection 
$\pi:X_{n,m}^1\to X_{n,m}$
\begin{equation*}
\begin{pmatrix}
a_{01}&\cdots&a_{0m}&x_0\\
a_{11}&\cdots&a_{1m}&x_1\\
\vdots &\cdots&\vdots &\vdots\\
a_{n1}&\cdots&a_{nm}&x_n\end{pmatrix}
\mapsto 
\begin{pmatrix}
a_{01}&\cdots&a_{0m}\\
a_{11}&\cdots&a_{1m}\\
\vdots &\cdots&\vdots \\
a_{n1}&\cdots&a_{nm}\end{pmatrix},
\end{equation*}
where each fiber is $\mpp^n$
and $X_{n,m}$ is referred to as the parameter space.

Let $A \in X_{n,m}$ and let $X_A$
be the rational variety obtained by
successive blow-ups at the points 
$P_i=P_i(A)=(a_{0i}:\cdots:a_{ni})$ $(i=1,2,\dots,m)$
from $\mpp^n$. We denote the family of rational projective varieties $X_A$ ($A\in X_{n,m}$) 
by $\wt{X}_{n,m}^1$, which also has the natural fibration  
$\wt{\pi}:\wt{X}_{n,m}^1 \to X_{n,m}$.

Let $E=E(A)$ be the divisor class on $X_A$ of 
the total transform of a hyper-plane in $\mpp^n$ and
let $E_i=E(A)$ be the exceptional divisor class generated by blow-up 
at the point $P_i$.
The group of divisor classes of $X_A$: $\pic(X_A) \simeq 
H^1(X_A,{\cal O}^{\times})
\simeq H^2(X_A,\mz)$ (the second equivalence comes from the fact that 
$X_A$ is a rational projective variety), is described as the lattice
\begin{equation}
\pic(X_A) = \mz E\oplus  \mz E_1 \oplus \mz E_2 \oplus \cdots \oplus 
\mz E_m.
\label{pic}
\end{equation}
Notice that this cohomology group is independent of $A$, while $X_A$ is
not isomorphic to $X_{A'}$ for $A' (\neq A) \in X_{n,m}$ in general. 

Let $e \in H_2(X_A,\mz)$ be the class of a generic line in $\mpp^n$
and let $e_i$ be the class of a generic line in the exceptional 
divisor of the blow-up at the point $P_i$.
Then, $e, e_1, e_2, \dots, e_m$ consist a basis of 
$H_2(X_A,\mz)\simeq (H^2(X_A,\mz))^*$ (the Poincar\'e duality) 
and the intersection numbers are given by 
$$\langle E,e \rangle =1, \quad
	\langle E,e_j \rangle =0, \quad
	\langle E_i, e \rangle =0, \quad
	\langle E_i,e_j \rangle = -\delta_{ij}.$$

Following Dolgachev-Ortland \cite{do}, 
we take the root basis $\{ \al_0,\dots, \al_{m-1}\}
\subset H^2(X_A,\mz)$
and the co-root basis $\{ \chal_0,\dots,\chal_{m-1}\} \subset H_2(X_A,\mz)$ as%
\ben &&\ba{ll}
\al_0= E - E_1-E_2-\cdots -E_{n+1}, & \al_i=E_i-E_{i+1} \ \ (i>0)\\
\chal_0= (n-1)e - e_1-e_2-\cdots -e_{n+1}, & \chal_i=e_i-e_{i+1} \ \ (i>0),
\ea \een
then, $\langle \al_i,\chal_i \rangle =-2$ holds for any $i$ and
these root bases define the Dynkin diagram of type $T_{2,n+1,m-n-1}$
by assigning a root $\al_i$ to every vertex $\al_i$  and
connecting two distinct vertices $\al_i$ and $\al_j$ 
if $\langle \al_i,\chal_j\rangle = 1$
(in our case $\langle \al_i,\chal_j\rangle = 0 \mbox{ or }1$ for $i\neq j$) 
(Fig.~\ref{dynkin1}).

Let us define the root lattice $Q=Q(n,m) \subset H^2(X_A,\mz)$ and 
the co-root lattice $Q^{\vee}=Q^{\vee}(n,m)\subset H_2(X_A,\mz)$ as
$Q=\mz \al_0 \oplus  \mz \al_1  \oplus \cdots \oplus \mz \al_{m-1} $
and $Q^{\vee}=\mz \chal_0 \oplus  \mz \chal_1  \oplus 
\cdots \oplus \mz \chal_{m-1} $ respectively.
For every $\al_i$ the formulae
\begin{equation}
	\begin{aligned}
{r_{\al_i}}_*(D)&={D+\langle D,\chal_i\rangle \al_i} \quad \mbox{for any }D\in Q \\
{r_{\al_i}}_*^{\vee}(d)&={d+ \langle \al_i,d \rangle \chal_i} \quad \mbox{for any } d \in Q^{\vee}
	\end{aligned}
	\label{simref}
\end{equation}
define linear involutions (called simple reflections) of the bi-lattice 
$(Q,Q^{\vee})$  
and they generate the Weyl group $W$ of type $T_{2,n+1,m-n-1}$,
which we denote by $W_*(n,m)$.

 
These simple reflections correspond to certain birational transformations
on the fiber space $\wt{\pi}: \wt{X}_{n,m}^1\to X_{n,m}$.
Let us define birational transformations 
$r_{i,j}$ $(1\leq i<j \leq m)$ and 
$r_{i_0,i_1,\dots, i_n}$ $(1\leq i_0<\cdots <i_n \leq m)$ 
on the fiber space as:\\
$r_{i,j}$ exchanges the points $P_i$ and $P_j$:
\begin{equation}
r_{i,j}: 
(\ba{c|c|c|c|c|c} \cdots&{\bf a}_i& \cdots& {\bf a}_j&\cdots&{\bf x} 
\ea) \mapsto
( \ba{c|c|c|c|c|c} \cdots& {\bf a}_j& \cdots& {\bf a}_i& \cdots &{\bf x}
\ea )\label{rij}
\end{equation}
and $r_{i_0,i_1,\dots, i_n}$ is the standard Cremona transformation
with respect to the points $P_{i_0}, P_{i_1}, \dots, P_{i_n}$, 
i.e. for example, $r_{1,2,\dots, n+1}$ is the composition of 
a projective transformation and the standard Cremona transformation
with respect to the origins $(0:\cdots:0:1:0\cdots:0)$ 
as
\begin{equation}
\begin{split}
r_{1,2,\dots, n+1}\colon
& (A~|~{\bf x})=
(\ba{c|c|c}A_{1,\dots, n+1}& A_{n+2,\dots, m}& {\bf x}\ea) \\
& \quad \mapsto
A_{1,\dots, n+1}^{-1} (\ba{c|c} A& {\bf x}\ea) 
=:\left(\ba{c|ccc|c}
 &       &\vdots  &        &\vdots \\
I_{n+1}&\cdots & a_{ij}''& \cdots &x_i''  \\
 &       &\vdots  &        &\vdots
\ea \right) \\
& \quad \mapsto 
\left(\ba{c|ccc|c}
 &       &\vdots  &        &\vdots \\
I_{n+1}&\cdots & {a_{ij}''}^{-1}& \cdots &{x_i''}^{-1} \\
 &       &\vdots  &        &\vdots
\ea \right), 
\end{split}
\label{r123}
\end{equation}
where $A_{j_1,j_2,\dots,j_k}$ denotes the $(n+1)\times k$ matrix
$({\bf a}_{j_1}~|~\cdots ~|~{\bf a}_{j_k})$
and $I_k$ denotes the $k \times k$ identity matrix.
(In section~5, ${a_{ij}''}^{-1}$ and ${x_i''}^{-1}$ are denoted as 
$a_{ij}'$ and $x_i'$, respectively.)

\medskip

Let $w$ denotes the reflection $r_{i,j}$ or $r_{i_0,i_1,\dots, i_n}$.
The reflection $w$ acts on the parameter space $X_{n,m}$ and 
preserves the fibration 
$\wt{\pi}:\wt{X}_{n,m}^1\to X_{n,m}.$ 
Recall that
$H^2(X_A,\mz)$ is independent of $A\in X_{n,m}$.
Hence, $w$ defines an action on this co-homology group.
Moreover, the induced birational map 
$w:X_A \dashrightarrow X_{w(A)}$ for generic $A\in X_{n,m}$ 
is a pseudo-isomorphism,
i.e. an isomorphism except sub-manifolds of co-dimension 2 or higher,
and the lines corresponding to the classes $e$ and $e_i$ can be chosen 
so that they do not meet the excepted part.
Since $H_2(X_A,\mz)$ is also independent of $A\in X_{n,m}$,
$w$ defines an action on this homology group and 
preserves the intersection form 
$\langle \cdot, \cdot \rangle : H^2(X_A,\mz)\times H_2(X_A,\mz)\to \mz$.

The birational maps $r_{i,i+1}$ and $r_{1,2,\dots,n+1}$ correspond 
to the simple reflections ${r_{\al_i}}_*$ $(1\leq i \leq m-1)$ and 
${r_{\al_0}}_*$ respectively. Indeed, their push-forward actions 
$H^2(X_A,\mz) \to H^2(X_{w(A)},\mz)$ 
$(w=r_{i,i+1}$ or $r_{1,2,\dots,n+1})$
and $H_2(X_A,\mz)\to H_2(X_{w(A)},\mz)$ 
are given by the formulae:
\begin{equation}
\begin{aligned}
{r_{i,i+1}}_*(D)&
=D+\langle D,\chal_i \rangle \al_i \\
{r_{i,i+1}}_*(d)&
=d+\langle \al_i,d \rangle \chal_i \\ 
{r_{1,2,\dots,n+1}}_*(D)&
=D+\langle D,\chal_0 \rangle \al_0 \\
{r_{1,2,\dots,n+1}}_*(d)&
=d+\langle \al_0,d \rangle \chal_0
\end{aligned}
\label{weylact}
\end{equation}
for any $D\in H^2(X_A,\mz)$ and any $d\in H_2(X_A,\mz)$.
Recall that root lattice $Q(n,m)$ and co-root lattice 
$Q^{\vee}(n,m)$ are subsets of
 $H^2(X_A,\mz)$ and $H_2(X_A,\mz)$ respectively and hence
the formulae (\ref{weylact}) are extensions of (\ref{simref}) onto 
these (co)-homology groups. \\

\definition
A hyper-surface in $X_A$ is called nodal  
if its class is $w_*(\al_0)$ for some $w_* \in W_*(n,m)$.
Let $N_{n,m}$ denote the set of $A \in X_{n,m}$ such that 
$X_A$ admits a nodal hyper-surface. 
We also write $N_{n,m}^1$ and $\wt{N}_{n,m}^1$ as
$\pi^{-1}(N_{n,m})$ and $\wt{\pi}^{-1}(N_{n,m})$ respectively.
\\

We define $W(n,m)$ as 
the group generated by  $r_{i,i+1}$ $(i=1,2,\dots,m-1)$ and 
$r_{1,2,\dots,n+1}$.

\begin{prop}[Coble, Dolgachev-Ortland] \label{cdo}
Let $m \geq n+2$. \\
{\rm i)} \ $W(n,m)\simeq W_*(n,m)$ holds. $W(n,m)$ acts on 
$\wt{X}_{n,m}^1 \setminus \wt{N}_{n,m}^1$.
 \\
{\rm ii)} \ Each element $w\in W(n,m)$ defines an action on 
$H^2(X_A,\mz)$ and $H_2(X_A,\mz)$, and 
preserves the intersection form 
$\langle \cdot , \cdot \rangle : H^2(X_A,\mz)\times H_2(X_A,\mz)\to \mz$.\\
{\rm iii)} \
The birational maps $r_{i,i+1}$ and $r_{1,2,\dots,n+1}$ correspond 
to the simple reflections $r_{\al_i}$ $(1\leq i \leq m-1)$ and 
$r_{\al_0}$ respectively.
For the reflection $r_{\al}=w\circ r_{\al_i} \circ w^{-1}$,
where $\al=w(\al_i)$ is a real root,
the formulae 
\begin{equation}
	\begin{aligned}
{r_{\al}}_*(D)&=D+\langle D,\chal \rangle \al\\
{r_{\al}}_*(d)&=d+\langle \al,d \rangle \chal
	\end{aligned} 
\label{weylact2}
\end{equation}
hold for any $D\in H^2(X_A,\mz)$ and any $d\in H_2(X_A,\mz)$.\\
{\rm iv)} \ Every $r_{i,j}$ or $r_{i_0,\dots,i_n}$ is an element of $W(n,m)$.
\end{prop}


\section{On the embedding of elliptic curve to $\mc \mpp^n$}

Let $T$ be an elliptic curve $\mc/(\mz+\mz \tau)$,
and let $\iota:T\to \mpp^n$ be an embedding.
We define the degree of $\iota$ as that of the pull-back 
of the line bundle ${\cal O}_{\mpp^n}(1)\simeq E$ by $\iota$.

\begin{prop}\label{prop310} {\rm (Algebraic version)}
Let $\iota$ and $\iota'$ be embeddings of $T$ to $\mpp^n$
of degree $n+1$ s.t. both $\iota(T)$ and $\iota'(T)$ are not
contained in any hyper-plane.
Then, there exists a translation $\sigma:T\to T$
and a projective linear transformation $G\in {\rm PGL}(n+1)$ s.t.
$G\circ \iota=\iota' \circ \sigma$ holds.
\end{prop}

$$\xymatrix{
\ar @{}[dr] |{\circlearrowleft}
\mpp^n \ar[r]^G & \mpp^n \\
T \ar[u]^\iota \ar[r]_\sigma  & T \ar[u]_{\iota'}
}$$

This proposition can be also stated as follows.

Let $f(u)$ be a holomorphic function on $\mc$.
We say that $f(u)$ has the quasi-periodicity if  
there exist constants $l_1,l_{\tau},c_1,c_{\tau}$ in $\mc$
s.t. the formulae
\begin{align*}
f(u+1)&=f(u)\exp\{2\pi \sqrt{-1}(l_1u+c_1)\}\\
f(u+\tau)&=f(u)\exp\{2\pi \sqrt{-1}(l_{\tau}u+c_{\tau})\}
\end{align*}
hold, and we refer such a function to as
a theta function for $T=\mc/(\mz+\mz \tau)$.

\begin{prop}{\rm (Analytic version)}\ 
Let $\iota$ be a holomorphic map 
$\iota:\mc \to \mpp^n: u \mapsto (f_0(u):f_1(u):\cdots:f_n(u))$ 
s.t.
{\rm (i)} $f_i$'s have the same quasi-periodicity, 
{\rm (ii)} $f_i$'s are linearly independent, 
and {\rm (iii)} each of them has $n+1$ zero points in the fundamental domain. 
Let
$\iota':\mc \to \mpp^n: u \mapsto (f_0'(u):f_1'(u):\cdots:f_n'(u))$ 
be also a holomorphic map s.t. $f_i$'s satisfy
{\rm (i)}, {\rm (ii)}, {\rm (iii)} with the same $\tau$
(the quasi-periodicity need not be the same with that of $f_i$'s).
Let $u_0,u_1,\dots,u_n$ be zeros of $f_0(u)$
and let $u_0',u_1',\dots,u_n'$ be zeros of $f_0'(u)$.
Then, there exists a projective linear transformation 
$G\in {\rm PGL}(n+1)$ s.t.
$\iota(u) G=\iota'(u+a)$ holds, where $a$ is determined by
\[
a=\frac{1}{n+1}(u_0'+u_1'+\cdots+u_n'-u_0-u_1-\cdots-u_n).
\] 
\end{prop}

\bp
The sum of zeros of $f_0'(u+a)$ is  
$(u_0'-a)+(u_1'-a)+\cdots+(u_n'-a)= u_0+u_1+\cdots+u_n$, then
 $f_0(u)$ and $f_0'(u+a)$ coincide up to a trivial theta function, i.e.
the quasi-periods of $f_0(u)$ and $f_0'(u+a)\exp(\al_2 u^2+\al_1 u)$ 
coincide for some $\al_1,\al_2\in \mc$.
 Thus, the quasi-periods of $f_i(u)$'s and 
$f_i'(u+a)\exp(\al_2u^2+\al_1u)$'s also coincide.
Hence, we can assume  $f_i$'s and $f_i'$'s have the same quasi-periodicity.

From the Riemann-Roch theorem
$$\dim H^0(T,{\cal O}(D))-\dim H^1(T,{\cal O}(D))=1-g+\deg D,$$ 
where $D$ is a divisor on $T$,
and the Serre duality $H^1(T,{\cal O}(D))^{*} \simeq H^0(T,\Omega(-D))$, 
we have
$$\dim H^0(T,{\cal O}(D))-\dim H^0(T,\Omega(-D))=\deg D.$$ 
Put $D=u_0+u_1+\cdots+u_n$ (summation of divisors), then we have
$H^0(T,\Omega(-D))=0$ and therefore $\dim H^0(T,{\cal O}(D))=n+1$.
Hence, meromorphic functions 
$f_i(u)/f_0(u)$ ($i=0,1,\dots,n$) on $T$ consist a basis
of $H^0(T,{\cal O}(D))$ and
$f_i'(u)/f_0(u)$ ($i=0,1,\dots,n$) also consist a basis
of $H^0(T,{\cal O}(D))$.
Hence, there exists a linear transformation $G$ on $\mc^{n+1}$
s.t. 
$$\left(\frac{f_0'(u)}{f_0(u)},\dots,\frac{f_n'(u)}{f_0(u)}\right)
=\left(\frac{f_0(u)}{f_0(u)},\dots,\frac{f_n(u)}{f_0(u)}\right)G $$
holds, and therefore,
$$\left(f_0'(u),\dots,f_n'(u)\right)
=\left(f_0(u),\dots,f_n(u)\right)G$$
holds. \ep


\section{Elliptic curves and birational 
representation of Weyl groups}

Let $X_A$ be the rational variety obtained by successive blow-ups of 
$\mpp^n$ at points $P_i$: $\rho_A:X_A\to \mpp^n$.   
Assume that there exist $\tau$, ${\rm Im} \, \tau>0$, and an embedding 
$\iota_A:T=\mc/(\mz+\mz \tau)\to \mpp^n$ s.t.
(i) the degree of $\iota_A$ is $n+1$;
(ii) all $P_i$ are on $\iota_A(T)$
(thus, $\iota_A(T)$ is not contained in any hyper-plane).
The embedding $\iota_A$ can be lifted to an embedding uniquely
$\tilde{\iota}_A :T \to X_A$ s.t. 
$\rho_A \circ \tilde{\iota}_A=\iota_A$ holds.
We denote by $\tilde{\iota}_A^*$ the pull-back 
$\tilde{\iota}_A^*:\pic(X_A)\to \pic(T)$.\\

\begin{prop}\label{prop410}
For $w\in W(n,m)$, $w$ induces an isomorphism from
$\tilde{\iota}_A(T)$ to $w\circ \tilde{\iota}_A(T)$.
\end{prop}

\bp
Since $P_{i_0},P_{i_1},\dots, P_{i_n}$ ($1\leq i_0 <i_1< \cdots <i_n\leq m$)
are not on any hyper-plane, the elliptic curve $\iota_A(T)\subset \mpp^n$ 
is not contained in a hyper-plane; therefore the generators 
$r_{1,2,\dots,n+1}$, $r_{i,i+1}$ act $\iota_A(T)$ birationally. 
Further, they can be extended to an isomorphism of elliptic curves.
Thus, by composition and lifting, the assertion follows.    
\ep

\begin{lemma}\label{lemma410}
The homology class of $\tilde{\iota}_A(T)$ in $H_2(X_A,\mz)$
is $(n+1)e-e_1-e_2-\cdots-e_m$.
\end{lemma}

\bp
The intersection numbers of $\tilde{\iota}_A(T)$ and the basis of 
$H^2(X_A,\mz)$ are
\begin{align*}
\langle E, \tilde{\iota}_A(T) \rangle &=n+1 \quad \mbox{(the degree of $\iota_A(T)$)}\\
\langle E_i, \tilde{\iota}_A(T) \rangle &=1 \quad (1\leq i \leq m).
\end{align*}
\ep\\

\begin{prop}\label{prop420}
The degree of $w\circ \iota_A:T\to \mpp^n$ is $n+1$.
\end{prop}

\bp
From (\ref{weylact}), we have
\begin{equation*}
\begin{split}
\deg(w\circ \iota_A(T))
&= {}_{H^2(X_{w(A)},\mz)}
\langle E, w\circ \tilde{\iota}_A(T)\rangle_{H_2(X_{w(A)},\mz)}\\
&= \langle E, w_*((n+1)e-e_1-e_2-\cdots-e_m)\rangle \\
&= \langle E, (n+1)e-e_1-e_2-\cdots-e_m)\rangle \\
&= n+1.
\end{split}
\end{equation*}
\ep\\

\noi {\it Notation}
For $u\in T$, we denote the divisor on $T$ corresponding to $u$ by $u$ again
and we denote its class by $[u]$. We denote the addition of
divisor classes $[u_1]$ and $[u_2]$ by $[u_1]+[u_2]$.\\

Let $Y_{n,m}$ denote a subset of $X_{n,m}$:
$$\{A\in X_{n,m}\setminus N_{n,m}~;~ 
\exists \tau, \mathrm{Im} \, \tau >0, \exists  \iota_A:T\to \mpp^n 
\mbox{: an embedding s.t.}\   (*)   \}$$
where $(*)$ is
(i) the degree of $\iota_A$ is $n+1$; (ii) all $P_i$ are on $\iota_A(T)$.
From Prop.~\ref{prop410} and Prop.~\ref{prop420}, the action of $W(n,m)$ 
can be restricted 
on the fiber space over $Y_{n,m}$.
For an embedding of an elliptic curve $\iota_A:T\to \mpp^n$ 
and $w\in W(n,m)$,
we write $w \circ \iota_A:T\to \mpp^n$ as $\iota_{w(A)}$.
It should be noted that $\iota_A$ and $\iota_{w(A)}$ are not
determined only by $A$ and $w(A)$, respectively,
e.g., both $\iota_A=(1,\wp(u),\wp'(u),\dots, \wp^{(n-1)}(u))$
and $\iota_A'=(1,\wp(-u),\wp'(-u),\dots, \wp^{(n-1)}(-u))$ satisfy
(i) and (ii). 
$$\xymatrix{
X_A \ar@{-->}[rr]^w  & &\ X_{w(A)}  \\
T \ar[u]^{\iota_A} \ar[urr]_{\iota_{w(A)}}
}$$

In this section, we consider the action of $W(n,m)$ on $Y_{n,m}$
via the orbit of an embedding of an elliptic curve $\iota_A:T\to \mpp^n$
for a parameter $A\in Y_{n,m}$.\\

We use the following notation ($**$):\\
$u$: a point on $T$;\\
$P_i$: a point in $\mpp^n$ determined by the 
$i$-th column of the parameter $A$;\\
$u_i$ ($1\leq i \leq m$): the point on $T$ s.t. $\iota_A(u_i)=P_i$;\\
$E$: the total transform of the class of a hyper-plane in $X_A$;\\  
$v$: a point in $T$ s.t. $(n+1)[v]=\tilde{\iota}_A^*(E)$
($\tilde{\iota}_A^*(E)$ is a line bundle on $T$ of degree
$n+1$, and therefore, there exists such a point $v\in T$.);\\
$A':=w(A)$;\\
$P_i'$: a point in $\mpp^n$ determined by the
$i$-th column of the parameter $A'$;\\
$u':=\iota_{A'}^{-1}\circ w\circ \iota_A(u)=u$.\\
$u_i'$ ($1\leq i \leq m$): the point on $T$ s.t. $\iota_A(u_i')=P_i'$;\\
$E'$: the total transform of the class of a hyper-plane in $X_{A'}$;\\
$v'$: a point in $T$ s.t. $(n+1)[v']=\tilde{\iota}_{A'}^*(E')$;\\

$$\xymatrix{
X_A \ar@{-->}[rr]^w  & & \ X_{A'}  \\
T \ar[u]^{\iota_A} \ar[urr]_{\iota_{A'}} &
} \quad
\xymatrix{
\pic(X_A) \ar[d]_{\iota_A^*} & &
\pic(X_{A'}) \ar[ll]_{w^*} \ar[dll]^{\iota_{A'}^*} \\
\pic(T) &}
$$

\begin{remark}
In the above diagrams, $P_i'\neq w(P_i)$ may occur, and therefore
$u_i'\neq \iota_{A'}^{-1} \circ w \circ \iota_A(u_i)$ also may occur. 
For example, we have $r_{ij}(P_i)=P_i$ and $P_i'=P_j$
($P_i'$ is the $i$-th column of $r_{ij}(A)$). 
\end{remark}

\begin{thm}\label{thm410}
Suppose that $w^*(E')$ and $w^*(E_i')$ are represented as 
$b_0^0 E+\sum_{j=1}^m b_0^j E_j$ and $b_i^0 E+\sum_{j=1}^m b_i^j E_j$,
respectively.
The points $u_i' \in T$ ($1\leq i \leq m$) and 
$u' \in T$ are given 
by the formulae:
\begin{align}
u_i'&= (n+1)b_i^0v+\sum_{j=1}^m b_i^j u_j \label{u_i'} \\
u' &= u  \label{u'}.
\end{align}
Moreover,
\begin{equation} 
(n+1)v'= (n+1)b_0^0v+\sum_{j=1}^m b_0^j u_j \label{v'}
\end{equation}
holds.
\end{thm}

\bp
From the definition of $\iota_{w(A)}$ (\ref{u'}) is trivial.
From the relation
$\iota_A^* \circ w^*= \iota_{A'}^*: 
\pic(X_{A'})\to \pic(T)$,
the divisor class
\begin{equation*}
\iota_A^* \circ w^*(E_i)
=\iota_A^*( b_i^0E+\sum_{j=1}^m b_i^j E_j)
=(n+1)b_i^0 [v] +\sum_{j=1}^m b_i^j [u_j]
\end{equation*}
coincides with the class
\begin{equation*}
\iota_{A'}^*(E_i)=[u_i'],
\end{equation*}
therefore, we have (\ref{u_i'}) by Abel's theorem.
Similarly, since the divisor class
\begin{equation*}
\iota_A^* \circ w^*(E)=\iota_A^*( b_0^0E+\sum_{j=1}^m b_0^j E_j)
=(n+1)b_0^0 [v] +\sum_{j=1}^m b_0^j [u_j]
\end{equation*}
coincides with the class
\begin{equation*}
\iota_{A'}^*(E)=(n+1)[v'],
\end{equation*}
we have (\ref{v'}).
\ep


\section{Representatives}

In order to investigate the actions of the Weyl group,
we should choose a suitable representative of the parameter $A$.
In this section, we consider realizations of $A$ and  $\iota_A$,
and study normalizations of $w(A)$ by $\mathrm{PGL}(n+1)$.

Let $A \in Y_{n,m}$, then, there 
exist $\tau$ (${\rm Im} \, \tau >0$) and an embedding $\iota_A: 
T=\mc/(\mz+\mz \tau) \to \mpp^n$ s.t.
(i) the degree of $\iota_A$ is $n+1$; (ii)  $\iota_A(T)$ contains $P_i$.
We use the notation ($**$) in the previous section.\\

\begin{thm}
Suppose that $w^*(E')$ is represented as $b_0^0 E+\sum_{j=1}^m b_0^j E_j$.
Then, there exists a translation $\sigma:T\to T:
\sigma(u)=u-s$
and a projective linear transformation $G\in {\rm PGL}(n+1)$ s.t.
\begin{equation}
\iota_A\circ \sigma =G\circ \iota_{w(A)} \label{th510}
\end{equation}
holds, where $s\in T$
satisfies
\begin{equation}
(n+1)s = (n+1)b_0^0v+\sum_{j=1}^m b_0^j u_j  -(n+1) v . \label{sigma}
\end{equation}
\end{thm}

\begin{remark}
By equality (\ref{th510}) the following diagram commutes
$$\xymatrix{
X_A \ar@{-->}[r]^w   & \ X_{w(A)} \ar[r]^G 
& X_{G \circ w(A)} & \hspace{-1cm}\\
T \ar[u]^{\iota_A} \ar[ur]_{\iota_{w(A)}} \ar[rr]_\sigma  & & 
T \ar[u]_{\iota_A}
}$$

\end{remark}

\bp
The relation (\ref{th510}) follows from Prop.\ref{prop310}.
Thus, we have
\begin{equation*}
\begin{split}
(n+1)[v']
&=\iota_{w(A)}^*(E')\\
&=(G^{-1}\circ \iota_A\circ \sigma)^*(E')\\
&=(\iota_A\circ \sigma)^*(E') \\
&=\sigma^*((n+1)[v])\\
&=(n+1)[v+s];
\end{split}
\end{equation*}
therefore, $(n+1)(v'-v)=(n+1)s \in T$ holds.
Further, Theorem~\ref{thm410} implies the equality (\ref{sigma}). \ep\\

In order to compute the action of $w$ on a general point in $\mpp^n$,
we have to determine $G\in \mathrm{PGL}(n+1)$ explicitly.  
For that purpose, it is sufficient to compute $n+2$ points in $\mpp^n$
in general position.

\noi i) The points $G(P_i')$ $(i=1,2,\dots,m)$ are calculated as
\begin{equation*}
\begin{split}
G(P_i')&= G\circ \iota_{w(A)}(u_i')\\
&= \iota_A \circ \sigma (u_i') \\
&= \iota_A (u_i'-s),
\end{split}
\end{equation*}
where $s$ and $u_i'$ are given by (\ref{sigma}) and (\ref{u_i'}), respectively.

\noi ii) The point $G \circ \iota_{w(A)}(u)$ 
$(u \neq u_i)$ is calculated as
\begin{equation*}
\begin{split}
G\circ \iota_{w(A)}(u')
&= \iota_A \circ \sigma (u') \\
&= \iota_A (u'-s)\\
&= \iota_A (u-s) .
\end{split}
\end{equation*}
The last equality follows from (\ref{u'}).


\subsection{Example 1}  
Let us calculate $G$ for the embedding
$\iota_A(u)={}^t(1,\wp(u),\wp'(u),\dots,\wp^{(n-1)}(u))$.
In this case, we can choose $v$ as $v=0$. 

\noi i)\ For $r_{1,2,\dots,n+1}$ (\ref{r123}).
From (\ref{th510}), the action is given by:
\begin{equation*}
\begin{split}
r_{1,2,\dots, n+1}\colon 
& (A~|~{\bf x})=(\iota_A(u_1),\iota_A(u_1),\dots,\iota_A(u_m)~|~{\bf x})\\
&\quad \mapright{A_{1,\dots, n+1}^{-1}} \ 
A_{1,\dots, n+1}^{-1} (A~|~{\bf x})
=:\left(\ba{c|ccc|c}
 &       &\vdots  &        &\vdots \\
I_{n+1}&\cdots & a_{ij}''& \cdots &x_i'' \\
 &       &\vdots  &        &\vdots
\ea \right)\\
 &\quad \mapright{{\rm SCT}} 
\left(\ba{c|ccc|c}
 &       &\vdots  &        &\vdots \\
I_{n+1}&\cdots & {a_{ij}''}^{-1}& \cdots &{x_i''}^{-1} \\
 &       &\vdots  &        &\vdots
\ea \right)=(A'~|~{\bf x}')\\
&\quad \mapright{G}~
(\ol{A}~|~\ol{{\bf x}})=
(\iota_A(u_1'+s),\iota_A(u_2'+s),\dots,\iota_A(u_m'+s)~|~\ol{{\bf x}}),
\end{split}
\end{equation*}
where 
$u_i'$ and $s$ are calculated by Theorem~\ref{thm410} 
and (\ref{sigma}) as
\begin{align*} u_i'&= 
\begin{cases}
\disp{u_i-\sum_{j=1}^{n+1}u_j }&(1\leq i \leq n+1)\\
u_i &(n+2\leq i\leq m)
\end{cases},\\ 
s &=-\frac{n-1}{n+1}\sum_{j=1}^{n+1}u_j.
\end{align*}
Thus, $G\in \mathrm{PGL}(n+1)$ is determined by 
$$G(I_{n+1},(\iota_A(0))')=
(\iota_A(u_1'+s),\iota_A(u_2'+s),\dots,\iota_A(u_{n+1}'+s),\iota_A(s)),$$
where $(\iota_A(0))'$ is
\begin{equation*}
\begin{split}
&(\iota_A(0))'=\\
&\, \left(
\left| 
  \ba{c}0\\\vdots\\0\\1\ea 
    \iota_A(u_2) \cdots \iota_A(u_{n+1}) \right|^{-1},
\left| 
  \iota_A(u_1)
  \ba{c}0\\\vdots\\0\\1\ea
    \cdots \iota_A(u_{n+1})
  \right|^{-1},\dots,
\left| \iota_A(u_1) \cdots \iota_A(u_n)
  \ba{c}0\\\vdots\\0\\1\ea 
  \right|^{-1}\right).
\end{split}
\end{equation*}
Here, $G$ can be decomposed as $G_2\circ G_1$,
where 
$$G_1(I_{n+1},(\iota_A(0))')=
\left(
\begin{array}{cc}
I_{n+1}&
\begin{array}{c}
1\\
\vdots\\
1
\end{array}
\end{array}
\right).$$

Thus, $G_1$ and $G_2$ are explicitly written as
\begin{equation*} 
\begin{split}
& G_1=\\
&{\rm diag}
\left(
\left| 
  \ba{c}0\\\vdots\\0\\1\ea 
    \iota_A(u_2) \cdots \iota_A(u_{n+1}) 
\right|,
\left| 
  \iota_A(u_1)
  \ba{c}0\\\vdots\\0\\1\ea
    \cdots \iota_A(u_{n+1})
\right|,\ldots,
\left| \iota_A(u_1) \cdots \iota_A(u_n)
  \ba{c}0\\\vdots\\0\\1\ea 
\right|
\right)
\end{split}
\end{equation*}
and $G_2^{-1}= ({\rm diag}(\ol{A}^{-1}\iota_A(s)))^{-1}\ol{A}^{-1}$;
thus, $G_2=\ol{A}\ {\rm diag}(\ol{A}^{-1}\iota_A(s))$.

\noi ii) For $r_{i,j}$ (\ref{rij}).
we have $u_i'=u_j,\ u_j'=u_i,\ u_k'=u_k\ (k\neq i,j)$, $s=0$, and 
$G$ is the identity.\\


\subsection{Example 2}  

Let us calculate $G$ for the embedding proposed by Kajiwara 
{\it et al.} \cite{kmnoy}:
$$\iota_A(u)=\left(\frac{[u-u_1-\ve]}{[u-u_1]}:\cdots: 
\frac{[u-u_{n+1}-\ve]}{[u-u_{n+1}]} \right),$$
where $[z]$ is a theta function whose
zero points are $\mz+ \mz \tau$ with the order 1,
and $\ve$ is an extra-parameter.
It should be noted that $A_{1,\dots,n+1}=
(\iota_A(u_1),\dots,\iota_A(u_{n+1}))=I_{n+1}$ holds.
In this case, we chose $G$ in a manner that 
$\iota_{\ol{A}}:=G\circ w\circ \iota_A$ 
is written in the form 
$$\iota_{\ol{A}}(u)=\left(\frac{[u-u_1'-\ol{\ve}]}{[u-u_1']}:\cdots: 
\frac{[u-u_{n+1}'-\ol{\ve}]}{[u-u_{n+1}']} \right).$$

From the diagram
$$\xymatrix{
X_A \ar@{-->}[r]^w   & \ X_{w(A)} \ar[r]^G 
& X_{G \circ w(A)} & \hspace{-1cm} \\
T \ar[u]^{\iota_A} \ar[ur]_{\iota_{w(A)}} \ar[rr]_\sigma  & & 
T \ar[u]_{\iota_{\ol{A}}}
}$$
we have
\begin{equation*}
\begin{split}
(n+1)[v']
&=\iota_{w(A)}^*(E') \\
&= (G^{-1}\circ \iota_{\ol{A}}\circ \sigma)^*(E') \\
&=(\iota_{\ol{A}}\circ \sigma)^*(E') \\
&=\sigma^*((n+1)[\ol{v}]) \\
&=(n+1)[\ol{v}+s],
\end{split}
\end{equation*}
and $v$, $v'$ and $\ol{v}$ are given by
\begin{align*}
(n+1)v&=\ve + \sum_{i=1}^{n+1}u_i\\
(n+1)\ol{v}&=\ol{\ve} + \sum_{i=1}^{n+1}u_i'\\
(n+1)v'&=(n+1)b_0^0 v + \sum_{j=1}^{m}b_0^j u_j.
\end{align*}

\noi i)\ For $r_{1,2,\dots,n+1}$. From theorem~\ref{thm410}, we have
$$
u_i'=
\begin{cases}
 u_i+\ve & (1\leq i \leq n+1)\\
 u_i& (n+2\leq i \leq m)
\end{cases}
,$$
$$(n+1)s=-\ve-\ol{\ve}.$$
Here, we can choose $\ol{\ve}$ as $\ol{\ve}=-\ve$, then we have $s=0$.

Further,
$G\in \mathrm{PGL}(n+1)$ is determined by 
$$G\circ w(A_{1,\dots,n+2})=G(I_{n+1},(\iota_A(u_{n+2}))')=
(I_{n+1},\iota_{\ol{A}}(u_{n+2}'))$$
and 
\begin{equation*} (\iota_A(u_{n+2}))'=
\left(\frac{[u_{n+2}-u_1]}{[u_{n+2}-u_1-\ve]}:\cdots: 
\frac{[u_{n+2}-u_{n+1}]}{[u_{n+2}-u_{n+1}-\ve]} \right).
\end{equation*}
Hence, $G$ is the identity.

\noi ii)\ For $r_{k,k+1}$, we have
\begin{equation*}
u_i'=
\begin{cases}
 u_i & (i\neq k,k+1)\\
 u_{k+1} & (i=k)\\
 u_k & (i=k+1)
\end{cases},
\end{equation*}
and
\begin{equation*}
(n+1)s=
\begin{cases}
 \ve-\ol{\ve} & (k\neq n+1)\\
 \ve-\ol{\ve}+u_{n+1}-u_{n+2} & (k=n+1)
\end{cases}.
\end{equation*}
Here, we can choose $\ol{\ve}$ as 
\begin{equation*}
\ol{\ve}=\begin{cases}
 \ve & (k\neq n+1)\\
 \ve+u_{n+1}-u_{n+2} & (k=n+1)
\end{cases} ,
\end{equation*}
then we have $s=0$.

For $r_{k,k+1}$ $(k\neq n+1)$, $G$ is the identity.
We calculate $G$ for $r_{n+1,n+2}$.
From
\begin{equation*}
G\left(
\ba{ccccc}
1&&&\frac{[u_{n+2}-u_1-\ve]}{[u_{n+2}-u_1]}&0\\
&\ddots&&\vdots&\vdots\\
&&1&\frac{[u_{n+2}-u_n-\ve]}{[u_{n+2}-u_n]}&0\\
&&&\frac{[u_{n+2}-u_{n+1}-\ve]}{[u_{n+2}-u_{n+1}]}&1
\ea \right)=
\left(
\ba{ccccc}
1&&&&\frac{[u_{n+2}-u_1-\ve]}{[u_{n+1}-u_1]}\\
&\ddots&&&\vdots\\
&&1&&\frac{[u_{n+2}-u_n-\ve]}{[u_{n+1}-u_n]}\\
&&&1&\frac{[-\ve]}{[u_{n+1}-u_{n+2}]}
\ea \right),
\end{equation*}
and by decomposing $G$ into $G_1=G_2\circ G_2$ as example 1,
we have
\begin{equation*} 
\begin{split}
G_1
&={\rm diag}
\left(
-\frac{[u_{n+2}-u_1]}{[u_{n+2}-u_1-\ve]},\dots,
-\frac{[u_{n+2}-u_n]}{[u_{n+2}-u_n-\ve]},1\right) \times \\
&\hspace{2cm} \times \left(
\ba{ccccc}
1&&&&
-\frac{[u_{n+2}-u_1-\ve][u_{n+2}-u_{n+1}]}{
[u_{n+2}-u_1][u_{n+2}-u_{n+1}-\ve]}\\
&\ddots&&&\vdots\\
&&1&&
-\frac{[u_{n+2}-u_n-\ve][u_{n+2}-u_{n+1}]}{
[u_{n+2}-u_n][u_{n+2}-u_{n+1}-\ve]}\\
&&&1&\frac{[u_{n+2}-u_{n+1}]}{[u_{n+2}-u_{n+1}-\ve]}
\ea \right)
\end{split}
\end{equation*}
and 
\begin{equation*}
G_2= 
\left(
\frac{[u_{n+2}-u_1-\ve]}{[u_{n+1}-u_1]},\dots,
\frac{[u_{n+2}-u_n-\ve]}{[u_{n+1}-u_n]},\frac{[-\ve]}{[u_{n+1}-u_{n+2}]}
\right).
\end{equation*}


\bigskip
\noindent {\it Acknowledgment.}
We would like to thank members of the Okamoto-Sakai seminar at Univ. of Tokyo for discussions and advice.
One of the authors (TT) appreciates the assistance from
Japan Society for the Promotion of Science.

\end{document}